\documentclass[journal]{IEEEtran}

\usepackage{amsmath, amssymb, epsfig}
\usepackage{fancyhdr}
\usepackage{bm}

\newlength{\algorithmwidth}
\newtheorem{theorem}{Theorem}[section]

\newtheorem{lemma}[theorem]{Lemma}

\newtheorem{definition}[theorem]{Definition}

\numberwithin{equation}{section}

\DeclareMathOperator*{\supp}{supp}

\DeclareMathOperator*{\argmin}{arg min}

\def \R {\mathbb{R}}

\def \e {\varepsilon}

\def \< {\langle}
\def \> {\rangle}
\def \^ {\widehat}

\def \supp {{\rm supp}}

\newcommand{\defby}{\overset{\mathrm{\scriptscriptstyle{def}}}{=}}

\newcommand{\bigO}{\mathrm{O}}

\begin{document}
\thispagestyle{empty}
\bibliographystyle{plain}
\title{Noisy Signal Recovery via Iterative Reweighted L1-Minimization}

\author{Deanna~Needell
\thanks{DN is with the Dept. of Mathematics, University of California, Davis, One Sheilds Ave., Davis CA 95616, USA e-mail: dneedell@math.ucdavis.edu.}
}

\date{April 2009}
\fancypagestyle{plain}
{
    \fancyhead{}
    \fancyfoot{}
}	
\thispagestyle{empty}
\maketitle
\thispagestyle{empty}
\begin{abstract}
Compressed sensing has shown that it is possible to reconstruct sparse high dimensional signals from few linear measurements. In many cases, the solution can be obtained by solving an $\ell_1$-minimization problem, and this method is accurate even in the presence of noise. Recent a modified version of this method, reweighted $\ell_1$-minimization, has been suggested. Although no provable results have yet been attained, empirical studies have suggested the reweighted version outperforms the standard method. Here we analyze the reweighted $\ell_1$-minimization method in the noisy case, and provide provable results showing an improvement in the error bound over the standard bounds.  
\end{abstract}

\section{Introduction}
Compressed sensing refers to the problem of realizing a sparse input $x$ using few linear measurements that possess some incoherence properties.  Its applications range from error correction to image processing.  Since the measurements are linear, the problem can be formulated as the recovery of a signal $x\in\R^d$ from its measurements $u=\Phi x$ where $\Phi$ is a $m \times d$ measurement matrix.  In the interesting case where $m \ll d$, it is clearly impossible to reconstruct any arbitrary signal, so we must restrict the domain to which the signals $x$ belong.  To that end, we consider \textit{sparse} signals, those with few non-zero coordinates relative to the actual dimension.  In particular, for $s\ll d$, we say that a signal $x\in\R^d$ is $s$-sparse if $x$ has $s$ or fewer non-zero coordinates:
$$
\|x\|_0 \defby |\{i : x_i \ne 0\}| \leq s.
$$
It is now well known that many signals such as real-world audio and video images or biological measurements are sparse either in this sense or with respect to a different basis.

Much work in the field of compressed sensing has led to promising reconstruction algorithms for these kinds of sparse signals. One solution to the recovery problem is simply to select the signal whose measurements are equal to those of $x$, with the smallest sparsity.  That is, one could solve the optimization problem
\begin{equation*}\tag{$L_0$}\label{L0}
\min_{\hat{x}\in\R^d} \|\hat{x}\|_0 \text{ subject to } \Phi \hat{x} = \Phi x.
\end{equation*}
This straightforward approach is quite accurate, and if the columns of the $m \times d$ matrix $\Phi$ are in general position, it can recover signals that are up to $m/2$-sparse.  The crucial drawback to the problem $(L_0)$ is of course that it is highly nonconvex and requires a search through the exponentially many column sets.  This clearly makes it of little use in practice.  A natural alternative then is to relax the problem $(L_0)$, and use a convex problem instead. One can then consider instead the $\ell_1$-minimization problem
\begin{equation*}\tag{$L_1$}\label{L1}
\min_{\hat{x}\in\R^d} \|\hat{x}\|_1 \text{ subject to } \Phi \hat{x} = \Phi x.
\end{equation*}
Here and throughout, $\|\cdot\|_1$ denotes the standard $\ell_1$-norm: $\|v\|_1 \defby \sum_{i=1}^d |v_i|$. 

Problem $(L_1)$ is convex, and can actually be reformulated as a linear program.  Due to recent work linear programming and smoothed analysis, it is now well known that it can be solved efficiently in practice~\cite{BV04:Convex,R06:Beyond}.  Notice that the solution to $(L_1)$ is the contact point where the smallest $\ell_1$-ball meets the subspace $x + \ker \Phi$.  The geometry of the octahedron lends itself well to sparsity due to its wedges at the lower dimensional subspaces. 

Indeed, Cand\`es and Tao prove that when the measurement matrix $\Phi$ satisfies a certain quantitative property, the solution to the problem $(L_1)$ will be the original sparse signal (\cite{CT05:Decoding}, see also~\cite{RV05:Geometric-Approach}).    This restricted isometry condition guarantees that every $m \times s$ submatrix of $\Phi$ approximately preserves norm: 

\begin{definition}
The measurement matrix $\Phi$ satisfies the restricted isometry condition with parameters $(s, \delta)$ if 
$$
(1-\delta)\|x\|_2^2 \leq \|\Phi x\|_2^2 \leq (1 + \delta)\|x\|_2^2,
$$
holds for all $s$-sparse vectors $x$. Here and throughout, $\|\cdot\|_2$ denotes the usual Euclidean norm: $\|x\|_2 \defby \big(\sum_{i=1}^d x_i^2 \big)^{1/2}$.
\end{definition}

It has been shown that many random measurement matrices satisfy the restricted isometry condition with small $\delta$ and $m$ nearly linear in the sparsity.  In ~\cite{MPJ06:Uniform} it is shown that measurement matrices whose entries are subgaussian satisfy the restricted isometry condition with parameters $(s, \delta)$ with high probability when 
$$
m = \bigO\big(\frac{s}{\delta^2}\log\frac{d}{\delta^2 s}\big).
$$
Note that this implies in particular that matrices whose entries are (normalized) random Gaussian or Bernoulli satisfy the restricted isometry condition with this number of measurements.  An alternative type of measurement matrix is a partial bounded orthogonal matrix.  One such example is obtained by selecting $m$ rows uniformly at random from the discrete Fourier matrix.  Rudelson and Vershynin show in~\cite{RV08:sparse} that these matrices satisfy the restricted isometry condition with parameters $(s, \delta)$ with
$$
m = \bigO\Big(\big(\frac{s\log d}{\epsilon^2}\big)\log\big(\frac{s\log d}{\epsilon^2}\big)\log^2 d \Big).
$$
Note that in both cases we need only $m \approx s\log d$ measurements.

Cand\`es and Tao showed that if the measurement matrix $\Phi$ satisfies the restricted isometry condition with parameters $(3s, 0.2)$, that the $s$-sparse signal $x$ is the unique solution to the problem $(L_1)$.  In~\cite{Can08:Restricted-Isometry}, Cand\`es sharpened these results to show success with parameters just $(2s, \sqrt{2}-1)$.

As is evident, these results provide strong guarantees for the $\ell_1$-minimization problem on sparse signals.  In practice, however, signals are rarely exactly sparse, and may often be corrupted by noise.  The problem then becomes to reconstruct an approximately sparse signal $x$ from its noisy measurements $u = \Phi x + e$, where $e$ is an error vector.  In this case, problem $(L_1)$ will clearly not suffice for recovery, since with noise the original signal may not even satisfy the constraint requirements.  However, the problem can simply be adapted to account for noise error:
\begin{equation*}\tag{$L_1'$}\label{L1p}
\min_{\hat{x}\in\R^d} \|\hat{x}\|_1 \text{ subject to } \|\Phi \hat{x} - u\|_2 \leq \varepsilon,
\end{equation*}
where $\varepsilon$ is a noise parameter with $\|e\|_2 \leq \varepsilon$.  Cand\`es, Romberg, and Tao showed that the solution to the problem $(L_1')$ is close in Euclidean norm to the original signal $x$~\cite{CRT06:Stable}.  Cand\`es improved these results to provide the following.

\begin{theorem}[$\ell_1$-minimization from~\cite{Can08:Restricted-Isometry}]\label{candes} Assume $\Phi$ has $\delta_{2s} < \sqrt{2}-1$. Let $x$ be an arbitrary signal with noisy measurements $\Phi x + e$, where $\|e\|_2 \leq \varepsilon$. Then the approximation $\hat{x}$ to $x$ from $\ell_1$-minimization satisfies
$$
\|x-\hat{x}\|_2 \leq C\varepsilon + C'\frac{\|x-x_s\|_1}{\sqrt{s}},
$$
where $C = \frac{2\alpha}{1-\rho}$, $C' = \frac{2(1+\rho)}{1-\rho}$, $\rho = \frac{\sqrt{2}\delta_{2s}}{1-\delta_{2s}}$ and $\alpha = \frac{2\sqrt{1+\delta_{2s}}}{\sqrt{1-\delta_{2s}}}$ .

\end{theorem}

The error bound provided here is optimal up to the constants, as the error $\varepsilon + \frac{\|x-x_s\|_1}{\sqrt{s}}$ can be viewed as the \textit{unrecoverable} energy due to the inherent noise.  See~\cite{NT08:Cosamp} for a detailed discussion of the unrecoverable energy.  Note also that in the case where $x$ is exactly sparse, but the measurements are noisy, the error bound $C\varepsilon$ is proportional to the norm of the error vector $e$.  Although these results provide very strong guarantees, recent work has been done on a variant of the $\ell_1$-minimization problem that seems to outperform the standard method.


\section{Reweighted $\ell_1$-minimization}\label{sec:rw}

As discussed above, the $\ell_1$-minimization problem $(L_1)$ is equivalent to the nonconvex problem $(L_0)$ when the measurement matrix $\Phi$ satisfies a certain condition.  The key difference between the two problems of course, is that $(L_1)$ depends on the magnitudes of the coefficients of a signal, whereas $(L_0)$ does not.  To reconcile this imbalance, a new weighted $\ell_1$-minimization algorithm was proposed by Cand\`es, Wakin, and Boyd~\cite{CWB08:Reweighted}.  This algorithm solves the following weighted version of $(L_1)$ at each iteration:

\begin{equation*}\tag{$WL_1$}\label{WL1}
\min_{\hat{x}\in\R^d} \sum_{i=1}^d \delta_i \hat{x}_i \text{ subject to } \Phi x = \Phi \hat{x}.
\end{equation*}

It is clear that in this formulation, large weights $\delta_i$ will encourage small coordinates of the solution vector, and small weights will encourage larger coordinates.  Indeed, suppose the $s$-sparse signal $x$ was known exactly, and that the weights were set as $\delta_i = \frac{1}{|x_i|}$.  Notice that in this case, the weights are infinite at all locations outside of the support of $x$.  This will force the coordinates of the solution vector $\hat{x}$ at these locations to be zero.  Thus if the signal $x$ is $s$-sparse with $s\leq m$, these weights would guarantee that $\hat{x} = x$.  Of course, these weights could not be chosen without knowing the actual signal $x$ itself.  However even if the weights are close to the actual signal, the geometry of the weighted $\ell_1$-ball becomes ``pinched'' toward the signal, decreasing the liklihood of an inaccurate solution.

Although the weights might not initially induce this geometry, one hopes that by solving the problem $(WL_1)$ at each iteration, the weights will get closer to the optimal values $\frac{1}{|x_i|}$, thereby improving the reconstruction of $x$.  Of course, one cannot actually have an infinite weight, so a stability parameter must also be used in the selection of the weight values.  The reweighted $\ell_1$-minimization algorithm can thus be described precisely as follows.

\bigskip

\textsc{Reweighted $\ell_1$-minimization}

\nopagebreak

\fbox{\parbox{\algorithmwidth}{
  \textsc{Input:} Measurement vector $u \in \R^m$, stability parameter $a$, noise parameter $\e$  
  
  \textsc{Output:} Reconstructed vector $\hat{x}$

  \begin{description}
    \item Initialize: Set the weights $\delta_i = 1$ for $i=1\ldots d$.\\
      Repeat the following until convergence or a fixed number of times:
    \item Approximate: Solve the reweighted $\ell_1$-minimization problem:
    $$ \hat{x} = \argmin_{\hat{x}\in\R^d} \sum_{i=1}^d \delta_i \hat{x}_i \text{ subject to } \|\Phi \hat{x} - u\|_2 \leq \e. $$
    \item Update the weights:
    $$ \delta_i = \frac{1}{|\hat{x}_i| + a}. $$
  \end{description}
 }}

    \bigskip


In~\cite{CWB08:Reweighted}, the reweighted $\ell_1$-minimization algorithm is discussed thoroughly, and experimental results are provided to show that it often outperforms the standard method.  However, no provable guarantees have yet been made for the algorithm's success.  Here we analyze the algorithm when the measurements and signals are corrupted with noise.  Since the reweighted method needs a weight vector that is somewhat close to the actual signal $x$, it is natural to consider the noisy case since the standard $\ell_1$-minimization method itself produces such a vector.  We are able to prove an error bound in this noisy case that improves upon the best known bound for the standard method.  We also provide numerical studies that show the bounds are improved in practice as well.



\section{Main Results}\label{sec:prfs}

The main theorem of this paper guarantees an error bound for the reconstruction using reweighted $\ell_1$-minimization that improves upon the best known bound of Theorem~\ref{candes} for the standard method.  For initial simplicity, we consider the case where the signal $x$ is exactly sparse, but the measurements $u$ are corrupted with noise.  Our main theorem, Theorem~\ref{theone} will imply results for the case where the signal $x$ is arbitrary.

\begin{theorem}[Reweighted $\ell_1$, Sparse Case]\label{theone}
Assume $\Phi$ satisfies the restricted isometry condition with parameters $(2s, \delta)$ where $\delta < \sqrt{2}-1$. Let $x$ be an $s$-sparse vector with noisy measurements $u = \Phi x + e$ where $\|e\|_2 \leq \varepsilon$. Assume the smallest nonzero coordinate $\mu$ of $x$ satisfies $\mu \geq \frac{4\alpha\varepsilon}{1-\rho}$.  Then the limiting approximation from reweighted $\ell_1$-minimization satisfies
$$
\|x-\hat{x}\|_2 \leq C''\e,
$$
where $C''= \frac{2\alpha}{1+\rho}$, $\rho = \frac{\sqrt{2}\delta}{1-\delta}$ and $\alpha = \frac{2\sqrt{1+\delta}}{1-\delta}$. 
\end{theorem}

\begin{remarks} 

{\bf 1.} We actually show that the reconstruction error satisfies
\begin{equation}\label{actualbnd}
\|x-\hat{x}\|_2 \leq \frac{2\alpha\varepsilon}{1 + \sqrt{1-\frac{4\alpha\varepsilon}{\mu}-\frac{4\alpha\varepsilon\rho}{\mu}}}.
\end{equation}
This bound is stronger than that given in Theorem~\ref{theone}, which is only equal to this bound when $\mu$ nears the value $\frac{4\alpha\varepsilon}{1-\rho}$.  However, the form in Theorem~\ref{theone} is much simpler and clearly shows the role of the parameter $\delta$ by the use of $\rho$.

{\bf 2.} For signals whose smallest non-zero coefficient $\mu$ does not satisfy the condition of the theorem, we may apply the theorem to those coefficients that do satisfy this requirement, and treat the others as noise.  See Theorem~\ref{ddnoise} below.

{\bf 3.} Although the bound in the theorem is the {\em limiting} bound, we provide a recursive relation~\eqref{Ek} in the proof which provides an exact error bound per iteration.  In Section~\ref{sec:num} we use dynamic programming to show that in many cases only a very small number of iterations are actually required to obtain the above error bound.
\end{remarks}


We now discuss the differences between Theorem~\ref{candes} and our new result Theorem~\ref{theone}. In the case where $\delta$ nears its limit of $\sqrt{2}-1$, the constant $\rho$ increases to $1$, and so the constant $C$ in Theorem~\ref{candes} is unbounded. However, the constant $C''$ in Theorem~\ref{theone} remains bounded even in this case.  In fact, as $\delta$ approaches $\sqrt{2}-1$, the constant $C''$ approaches just $4.66$.  The tradeoff of course, is in the requirement on $\mu$.  As $\delta$ gets closer to $\sqrt{2}-1$, the bound needed on $\mu$ requires the signal to have unbounded non-zero coordinates relative to the noise level $\e$.  However, to use this theorem efficiently, one would select the largest $\delta < \sqrt{2}-1$ that allows the requirement on $\mu$ to be satisfied, and then apply the theorem for this value of $\delta$.  Using this strategy, when the ratio $\frac{\mu}{\e} = 10$, for example, the error bound is just $3.85\e$.  

Theorem~\ref{theone} and a short calculation will imply the following result for \textit{arbitrary} signals $x$.

\begin{theorem}[Reweighted $\ell_1$]\label{ddnoise}
Assume $\Phi$ satisfies the restricted isometry condition with parameters $(2s, \sqrt{2}-1)$. Let $x$ be an arbitrary vector with noisy measurements $u = \Phi x + e$ where $\|e\|_2 \leq \varepsilon$. Assume the smallest nonzero coordinate $\mu$ of $x_s$ satisfies $\mu \geq \frac{4\alpha\varepsilon_0}{1-\rho},$ where $\varepsilon_0 = 
1.2(\|x-x_s\|_2 + \frac{1}{\sqrt{s}}\|x-x_s\|_1) + \varepsilon$.
Then the limiting approximation from reweighted $\ell_1$-minimization satisfies
$$
\|x-\hat{x}\|_2 \leq \frac{4.1\alpha}{1+\rho}\Big(\frac{\|x-x_{s/2}\|_1}{\sqrt{s}} + \varepsilon \Big),
$$
and
$$
\|x-\hat{x}\|_2 \leq \frac{2.4\alpha}{1+\rho}\Big(\|x-x_s\|_2 + \frac{\|x-x_s\|_1}{\sqrt{s}} + \varepsilon \Big),
$$
where $\rho$ and $\alpha$ are as in Theorem~\ref{theone}. 
\end{theorem}

  Again in the case where $\delta$ nears its bound of $\sqrt{2}-1$, both constants $C$ and $C'$ in Theorem~\ref{candes} approach infinity. However, in Theorem~\ref{ddnoise}, the constant remains bounded even in this case.  The same strategy discussed above for Theorem~\ref{theone} should also be used for Theorem~\ref{ddnoise}.  Next we begin proving Theorem~\ref{theone} and Theorem~\ref{ddnoise}.
  
  

\subsection{Proofs}

We will first utilize a lemma that bounds the $\ell_2$ norm of a small portion of the difference vector $x - \hat{x}$ by the $\ell_1$-norm of its remainder.  This lemma is proved in~\cite{Can08:Restricted-Isometry} and essentially in~\cite{CRT06:Stable} as part of the proofs of the main theorems of those papers.  

\begin{lemma}\label{tinyLem}
Set $h = \hat{x} - x$, and let $\alpha$, $\e$, and $\rho$ be as in Theorem~\ref{theone}. Let $T_0$ be the set of $s$ largest coefficients in magnitude of $x$ and $T_1$ be the $s$ largest coefficients of $h_{T_0^c}$. Then
\begin{equation}\label{lembeg}
\|h_{T_0\cup T_1}\|_2 \leq \alpha\e + \frac{\rho}{\sqrt{s}}\|h_{T_0^c}\|_1,
\end{equation}
and
\begin{equation}\label{ref1}
\|h_{(T_0\cup T_1)^c}\|_2 \leq \frac{1}{\sqrt{s}}\|h_{T_0^c}\|_1.
\end{equation}
\end{lemma}%

We will next require two lemmas that give results about a single iteration of reweighted $\ell_1$-minimization.

\begin{lemma}[Single reweighted $\ell_1$-minimization]\label{doublenoise}
Assume $\Phi$ satisfies the restricted isometry condition with parameters $(2s, \sqrt{2}-1)$. Let $x$ be an arbitrary vector with noisy measurements $u = \Phi x + e$ where $\|e\|_2 \leq \varepsilon$. Let $w$ be a vector such that $\|w-x\|_\infty \leq A$ for some constant $A$. Denote by $x_s$ the vector consisting of the $s$ (where $s \leq |\supp (x)|$) largest coefficients of $x$ in absolute value. Let $\mu$ be the smallest coordinate of $x_s$ in absolute value, and set $b = \|x-x_s\|_\infty$. Then when $\mu \geq A$ and $\rho C_1 < 1$, the approximation from reweighted $\ell_1$-minimization using weights $\delta_i = 1 / (w_i + a)$ satisfies
$$
\|x-\hat{x}\|_2 \leq
D_1\varepsilon + D_2\frac{\|x-x_s\|_1}{a},
$$
where $D_1=\frac{(1+C_1)\alpha}{1-\rho C_1}$, $D_2 = C_2+\frac{(1+C_1)\rho C_2}{1-\rho C_1}$, $C_1 = \frac{A+a+b}{\mu-A+a}$, $C_2 = \frac{2(A+a+b)}{\sqrt{s}}$, and $\rho$ and $\alpha$ are as in Theorem~\ref{theone}.
\end{lemma}

\begin{proof} Now we begin the proof of Lemma~\ref{doublenoise}.

Set $h$ and $T_j$ for $j\geq 0$ as in Lemma~\ref{tinyLem}.  For simplicity, denote by $\|\cdot\|_w$ the weighted $\ell_1$-norm: 
$$
\|z\|_w \defby \sum_{i=1}^d \frac{1}{|w_i| + a}z_i.
$$
Since $\hat{x} = x+h$ is the minimizer of~\eqref{WL1}, we have 
\begin{align*}
\|x\|_w &\geq \|x+h\|_w = \|(x+h)_{T_0}\|_w + \|(x+h)_{T_0^c}\|_w \\
&\geq \|x_{T_0}\|_w - \|h_{T_0}\|_w + \|h_{T_0^c}\|_w - \|x_{T_0^c}\|_w.
\end{align*}
This yields 
$$
\|h_{T_0^c}\|_w \leq \|h_{T_0}\|_w + 2\|x_{T_0^c}\|_w.
$$
Next we relate the reweighted norm to the usual $\ell_1$-norm. We first have
$$
\|h_{T_0^c}\|_w \geq \frac{\|h_{T_0^c}\|_1}{A+a+b},
$$
by definition of the reweighted norm as well as the values of $A$, $a$, and $b$.
Similarly we have
$$
\|h_{T_0}\|_w  \leq \frac{\|h_{T_0}\|_1}{\mu - A + a}.
$$
Combining the above three inequalities, we have
\begin{align}\label{ref2}
\|h_{T_0^c}\|_1 &\leq (A+a+b)\|h_{T_0^c}\|_w \\
&\leq (A+a+b)(\|h_{T_0}\|_w + 2\|x_{T_0^c}\|_w)\notag\\
&\leq \frac{A+a+b}{\mu-A+a}\|h_{T_0}\|_1 + 2(A+a+b)\|x_{T_0^c}\|_w.
\end{align}
Using~\eqref{ref1} and~\eqref{ref2} along with the fact $\|h_{T_0}\|_1 \leq \sqrt{s}\|h_{T_0}\|_2$, we have
\begin{equation}\label{third}
\|h_{(T_0\cup T_1)^c}\|_2 \leq C_1\|h_{T_0}\|_2 + C_2\|x_{T_0^c}\|_w, 
\end{equation}
where $C_1 = \frac{A+a+b}{\mu-A+a}$ and $C_2 = \frac{2(A+a+b)}{\sqrt{s}}$.
By~\eqref{lembeg} of Lemma~\ref{tinyLem}, we have
$$
\|h_{T_0\cup T_1}\|_2 \leq \alpha\varepsilon + \frac{\rho}{\sqrt{s}}\|h_{T_0^c}\|_1,
$$
where $\rho = \frac{\sqrt{2}\delta_{2s}}{1-\delta_{2s}}$ and $\alpha = \frac{2\sqrt{1+\delta_{2s}}}{\sqrt{1-\delta_{2s}}}$.
Thus by~\eqref{ref2}, we have 
\begin{align*}
\|h_{T_0\cup T_1}\|_2 &\leq \alpha\varepsilon + \frac{\rho}{\sqrt{s}}(C_1\|h_{T_0}\|_1 + 2(A+a+b)\|x_{T_0^c}\|_w) \\
&= \alpha\varepsilon + \rho C_1\|h_{T_0\cup T_1}\|_2 + \rho C_2\|x_{T_0^c}\|_w.
\end{align*}
Therefore, 
\begin{equation}\label{need}
\|h_{T_0\cup T_1}\|_2 \leq (1 - \rho C_1)^{-1}(\alpha\varepsilon + \rho C_2\|x_{T_0^c}\|_w).
\end{equation}
Finally by~\eqref{third} and~\eqref{need},
\begin{align*}
\|h\|_2 &\leq \|h_{T_0\cup T_1}\|_2 + \|h_{(T_0\cup T_1)^c}\|_2 \\
&\leq(1+C_1)\|h_{T_0\cup T_1}\|_2 + C_2\|x_{T_0^c}\|_w \\
&\leq (1+C_1)\Big(\frac{\alpha\varepsilon + \rho C_2\|x_{T_0^c}\|_w}{1 - \rho C_1}\Big) + C_2\|x_{T_0^c}\|_w.
\end{align*}
Applying the inequality $\|x_{T_0^c}\|_w \leq (1/a)\|x_{T_0^c}\|_1$ and simplifying completes the claim.

\end{proof}

Applying Lemma~\ref{doublenoise} to the case where $x-x_s = 0$ and $b=0$ yields the following. 

\begin{lemma}[Single reweighted $\ell_1$-minimization, Sparse Case]\label{noisymeas}
Assume $\Phi$ satisfies the restricted isometry condition with parameters $(2s, \sqrt{2}-1)$. Let $x$ be an $s$-sparse vector with noisy measurements $u = \Phi x + e$ where $\|e\|_2 \leq \varepsilon$. Let $w$ be a vector such that $\|w-x\|_\infty \leq A$ for some constant $A$. Let $\mu$ be the smallest non-zero coordinate of $x$ in absolute value. Then when $\mu \geq A$, the approximation from reweighted $\ell_1$-minimization using weights $\delta_i = 1 / (w_i + a)$ satisfies
$$
\|x-\hat{x}\|_2 \leq D_1\varepsilon. 
$$
Here $D_1 = \frac{(1+C_1)\alpha}{1-\rho C_1}$, $C_1 = \frac{A+a}{\mu-A+a}$, and $\alpha$ and $\rho$ are as in Theorem~\ref{theone}.

\end{lemma}

Now we begin the proof of Theorem~\ref{theone}.
\begin{proof} The proof proceeds as follows. First, we use the error bound in Theorem~\ref{candes} as the initial error, and then apply Lemma~\ref{noisymeas} repeatedly. We show that the error decreases at each iteration, and then deduce its limiting bound using the recursive relation. To this end, let $E(k)$ for $k=1,\ldots$, be the error bound on $\|x-\hat{x}_k\|_2$ where $\hat{x}_k$ is the reconstructed signal at the $k^{th}$ iteration. Then by Theorem~\ref{candes} and Lemma~\ref{noisymeas}, we have the recursive definition
\begin{equation}\label{Ek}
E(1) = \frac{2\alpha}{1-\rho}\varepsilon, \quad E(k+1) = \frac{(1+\frac{E(k)}{\mu-E(k)})\alpha}{1-\rho \frac{E(k)}{\mu-E(k)}}\varepsilon. 
\end{equation}
Here we have taken $a\rightarrow 0$ iteratively (or if $a$ remains fixed, a small constant $O(a)$ will be added to the error). First, we show that the base case holds, $E(1) \leq E(2)$. Since $\mu \geq \frac{4\alpha\varepsilon}{1-\rho}$, we have that
$$
\frac{E(1)}{\mu-E(1)} = \frac{\frac{2\alpha\varepsilon}{1-\rho}}{\mu-\frac{2\alpha\varepsilon}{1-\rho}} \leq 1.
$$
Therefore we have
$$
E(2) = \frac{(1+\frac{E(1)}{\mu-E(1)})\alpha}{1-\rho \frac{E(1)}{\mu-E(1)}}\varepsilon \leq \frac{2\alpha}{1-\rho}\varepsilon = E(1).
$$
Next we show the inductive step, that $E(k+1)\leq E(k)$ assuming the inequality holds for all previous $k$. Indeed, if $E(k) \leq E(k-1)$, then we have
$$
E(k+1) = \frac{(1+\frac{E(k)}{\mu-E(k)})\alpha}{1-\rho \frac{E(k)}{\mu-E(k)}}\varepsilon \leq \frac{(1+\frac{E(k-1)}{\mu-E(k-1)})\alpha}{1-\rho \frac{E(k-1)}{\mu-E(k-1)}}\varepsilon = E(k).
$$  
Since $\mu \geq \frac{4\alpha\varepsilon}{1-\rho}$ and $\rho \leq 1$ we have that $\mu - E(k)\geq 0$ and $\rho\frac{E(k)}{\mu-E(k)}\leq 1$, so $E(k)$ is also bounded below by zero. Thus E(k) is a bounded decreasing sequence, so it must converge. Call its limit $L$. By the recursive definition of $E(k)$, we must have
$$
L = \frac{(1+\frac{L}{\mu-L})\alpha}{1-\rho \frac{L}{\mu-L}}\varepsilon.
$$
Solving this equation yields
$$
L=\frac{\mu-\sqrt{\mu^2-4\mu\alpha\varepsilon-4\mu\alpha\varepsilon\rho}}{2(1+\rho)},
$$
where we choose the solution with the minus since $E(k)$ is decreasing and $E(1) < \mu/2$ (note also that $L=0$ when $\e=0$). Multiplying by the conjugate and simplifying yields
\begin{align*}
L &= \frac{4\mu\alpha\varepsilon+4\mu\alpha\varepsilon\rho}{2(1+\rho)(\mu+\sqrt{\mu^2-4\mu\alpha\varepsilon-4\mu\alpha\varepsilon\rho})} \\
&= \frac{2\alpha\varepsilon}{1 + \sqrt{1-\frac{4\alpha\varepsilon}{\mu}-\frac{4\alpha\varepsilon\rho}{\mu}}}.
\end{align*}
Then again since $\mu \geq \frac{4\alpha\varepsilon}{1-\rho}$, we have
$$
L \leq \frac{2\alpha\varepsilon}{1+\rho}.
$$ 
\end{proof}

\begin{proof} Now we begin the proof of Theorem~\ref{ddnoise}.
By Lemma 6.1 of~\cite{NT08:Cosamp} and Lemma 7 of~\cite{GSTV07:HHS}, we can rewrite $\Phi x + e$ as $\Phi x_s + \widetilde{e}$ where 
\begin{align*}
\|\widetilde{e}\|_2 &\leq 1.2(\|x-x_s\|_2 + \frac{1}{\sqrt{s}}\|x-x_s\|_1) + \|e\|_2 \\
&\leq 2.04\Big(\frac{\|x-x_{s/2}\|_1}{\sqrt{s}}\Big) + \|e\|_2.
\end{align*}
This combined with Theorem~\ref{theone} completes the claim.
\end{proof}

\section{Numerical Results and Convergence}\label{sec:num}

\begin{figure}[ht] 
  \includegraphics[scale=0.4]{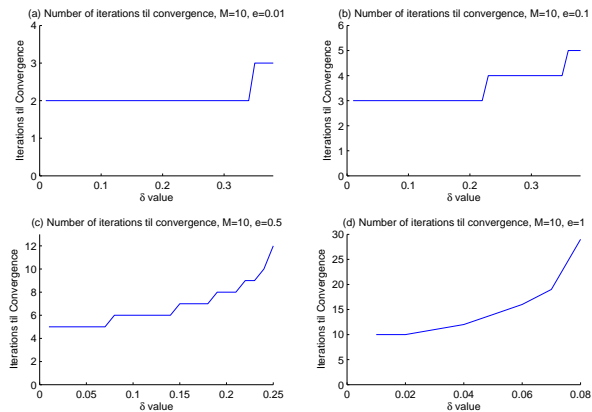}
  \caption{The number of iterations required for the theoretical error bounds~eqref{Ek} to reach the theoretical error threshold~\eqref{actualbnd} when (a) $\mu = 10$, $\e=0.01$, (b) $\mu = 10$, $\e=0.1$, (c) $\mu = 10$, $\e=0.5$, (d) $\mu = 10$, $\e=1.0$. }\label{fig:its}
\end{figure}

Our main theorems prove bounds on the reconstruction error limit.  However, as is the case with many recursive relations, convergence to this threshold is often quite fast.  To show this, we use dynamic programming to compute the theoretical error bound $E(k)$ given by~\eqref{Ek} and test its convergence rate to the threshold given by~eqref{actualbnd}. Since the ratio between $\mu$ and $\e$ is important, we fix $\mu = 10$ and test the convergence for various values of $\e$ and $\delta$.  The results are displayed in Figure~\ref{fig:its}.  We observe that in each case, as $\delta$ increases we require slightly more iterations.  This is not surprising since higher $\delta$ means a lower bound.  We also confirm that less iterations are required when the ratio $\mu/\e$ is smaller. 

\begin{figure}[h!]
\begin{center}
$\begin{array}{c@{\hspace{.1in}}c}
\includegraphics[scale=0.37]{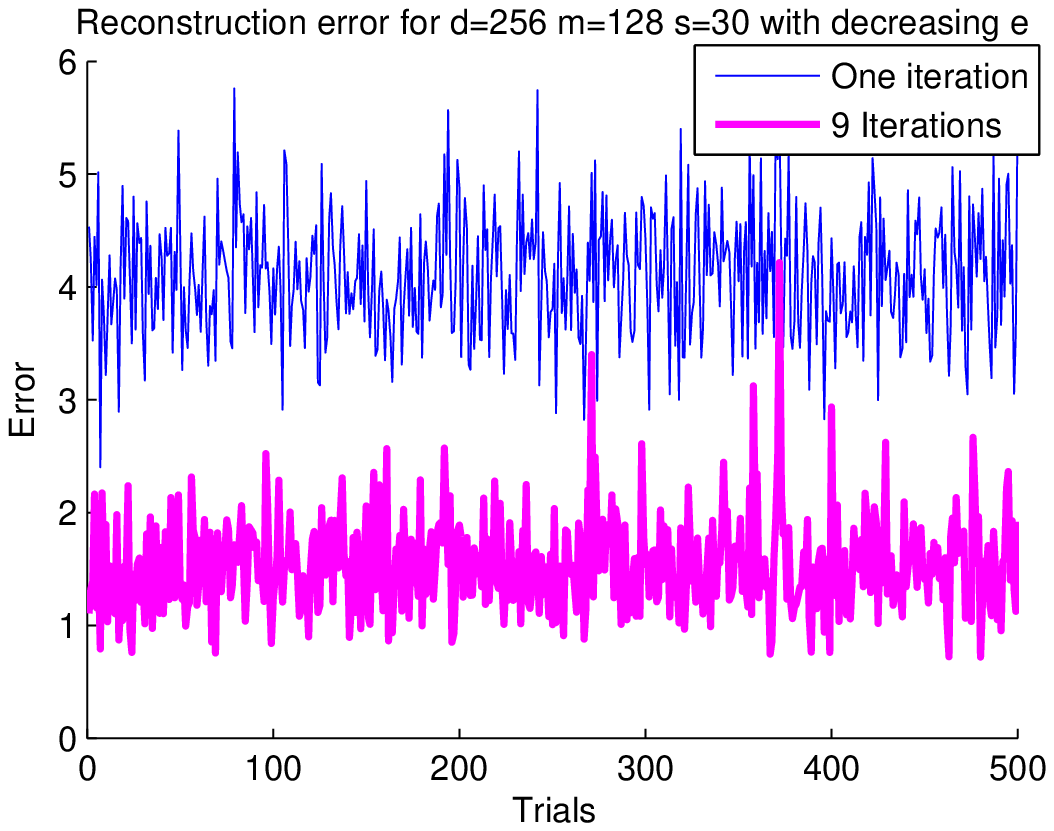}  &  
\includegraphics[scale=0.37]{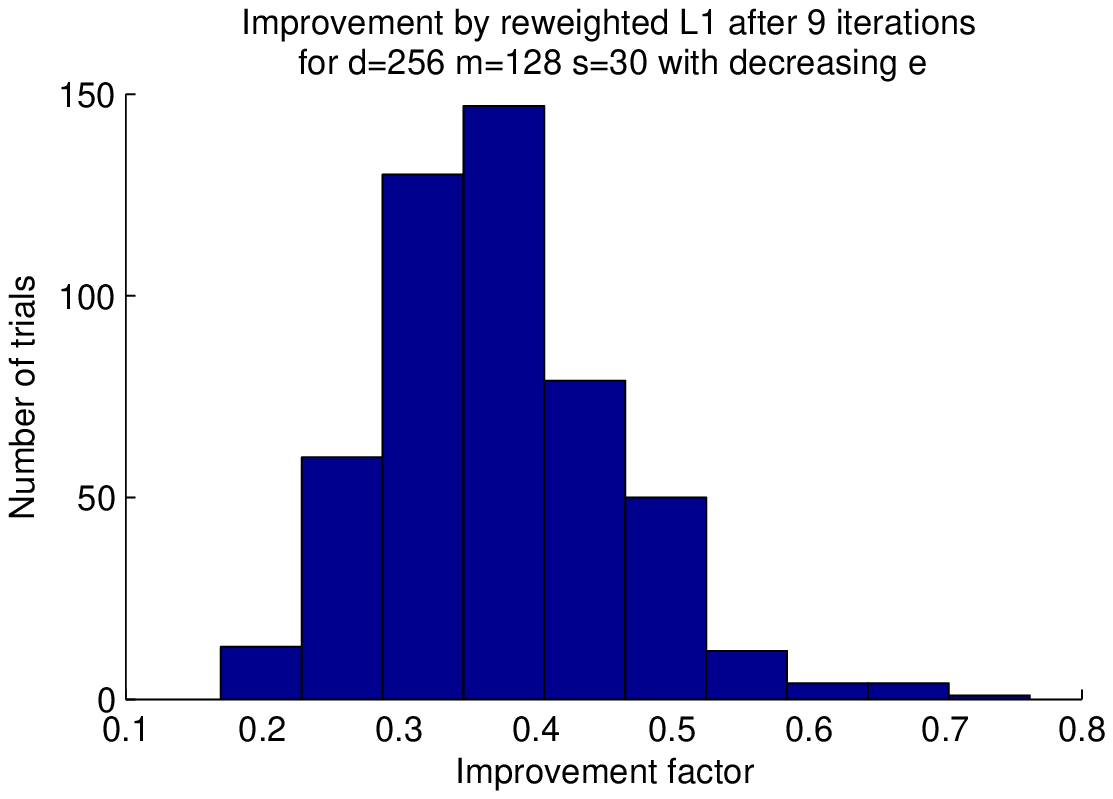}   \\
\end{array}$
\end{center}
\caption{Improvements in the $\ell_2$ reconstruction error using reweighted $\ell_1$-minimization versus standard $\ell_1$-minimization for Gaussian signals. Error plot (left) and histogram of improvement factors $\|x - \hat{x}\|_2 / \|x - x^*\|_2$ (right).}
\label{fig:dece}
\end{figure}

\begin{figure}[h!]
\begin{center}
$\begin{array}{c@{\hspace{.1in}}c}
\includegraphics[scale=0.37]{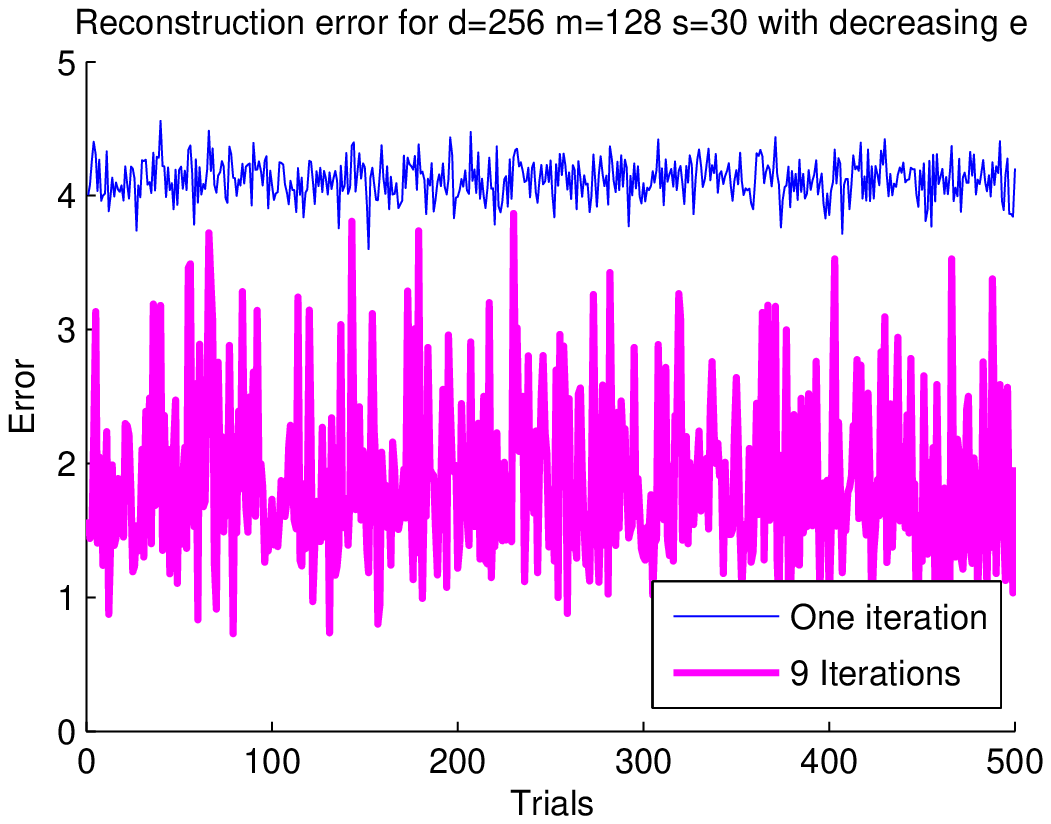}  &  
\includegraphics[scale=0.37]{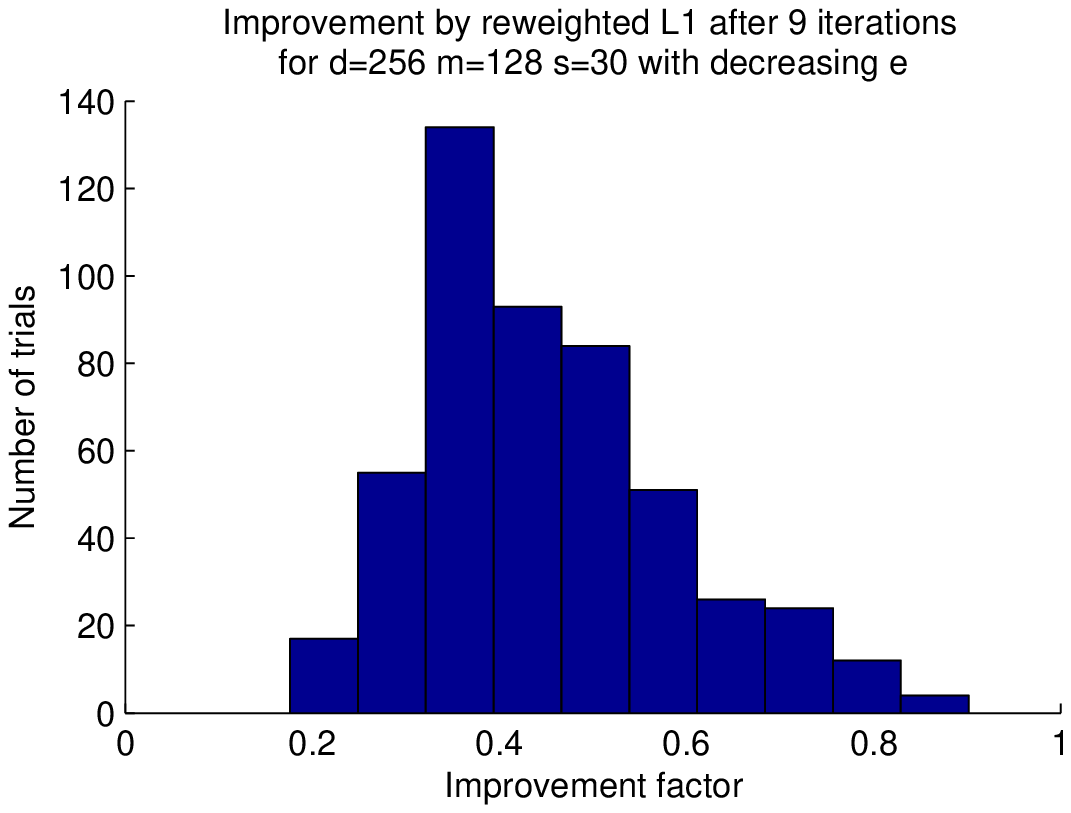}   \\
\end{array}$
\end{center}
\caption{Improvements in the $\ell_2$ reconstruction error using reweighted $\ell_1$-minimization versus standard $\ell_1$-minimization for Bernoulli signals. Error plot (left) and histogram of improvement factors $\|x - \hat{x}\|_2 / \|x - x^*\|_2$ (right).}
\label{fig:deceBern}
\end{figure}

Next we examine some numerical experiments conducted to test the actual error with reweighted $\ell_1$-minimization versus the standard $\ell_1$ method.  In these experiments we consider signals of dimension $d=256$ with $s=30$ non-zero entries.  We use a $128 \times 256$ measurement matrix $\Phi$ consisting of Gaussian entries.  We note that we found similar results when the measurement matrix $\Phi$ consisted of symmetric Bernoulli entries.  
For each trial in our experiments we construct an $s$-sparse signal $x$ with support chosen uniformly at random and entries from either the Gaussian distribution or the symmetric Bernoulli distribution, all independent of the matrix $\Phi$. We then construct the normalized Gaussian noise vector $e$, and run the reweighted $\ell_1$-algorithm using $\e$ such that $\e^2 = \sigma^2(m + 2\sqrt{2m})$ where $\sigma^2$ is the variance of the normalized error vectors. This value is likely to provide a good upper bound on the noise norm (see e.g. \cite{CRT06:Stable}, \cite{CWB08:Reweighted}). We set $a = k/1000$ in the $k^{th}$ iteration.  We run 500 trials for each parameter selection and signal type.  We found similar results for non-sparse signals, which is not surprising since we can treat the signal error as measurement error after applying the measurement matrix (see the proof of Theorem~\ref{ddnoise}). Figures~\ref{fig:dece} and~\ref{fig:deceBern} display the results of the experiments and demonstrate large improvements in the error of the reweighted reconstruction $\hat{x}$ compared to the reconstruction $x^*$ from the standard method.  

  
\bibliography{rw}

\end{document}